\documentclass[12pt]{article}
\usepackage{amsfonts}
\usepackage{amsmath}
\usepackage{amssymb}
\usepackage{graphicx}
\usepackage{indentfirst}
\usepackage{bm}
\usepackage{caption}
\usepackage[T1]{fontenc}
\usepackage{geometry}
\makeatletter

\newcommand{\Rmnum}[1]{\expandafter\@slowromancap\romannumeral #1@}
\makeatother

\begin{document}
\title{Improved mutual coherence of some random
overcomplete dictionaries for sparse repsentation\thanks{This work was supported by the National Natural Science Foundation of China under Grant 11131006. This work was supported in part by the National Basic Research Program of China under Grant 2013CB329404.}}
\author{%
\small{\textsf{Yingtong Chen}}\\
\texttt{\small{chenyingtong@stu.xjtu.edu.cn}}\\
\small{School of Mathematics and Statistics, Xi'an Jiaotong University,}\\
\small{No.28, Xianning West Road, Xi'an, Shaanxi, 710049, P.R. China}\\
\small{Beijing Center for Mathematics and Information}\\ \small{Interdisciplinary Sciences (BCMIIS)} \\
\small{\textsf{Jigen Peng}} \\
\texttt{\small{jgpeng@mail.xjtu.edu.cn}} \\
\small{School of Mathematics and Statistics, Xi'an Jiaotong University,}\\
\small{No.28, Xianning West Road, Xi'an, Shaanxi, 710049, P.R. China}\\
\small{Beijing Center for Mathematics and Information}\\
\small{Interdisciplinary Sciences (BCMIIS)}
}
\date{}
\maketitle

\begin{abstract}
The letter presents a method for the reduction in the mutual coherence of an overcomplete Gaussian or Bernoulli random matrix, which is fairly small due to the lower bound given here on the probability of the event that the aforesaid mutual coherence is less than any given number in (0, 1). The mutual coherence of the matrix that belongs to a set which contains the two types of matrices with high probability can be reduced by a similar method but a subset that has Lebesgue measure zero. The numerical results are provided to illustrate the reduction in the mutual coherence of an overcomplete Gaussian, Bernoulli or uniform random dictionary. The effect on the third type is better than a former result.
\end{abstract}
\newenvironment{keywords}{\par\noindent{\bf Keywords:}}{\par}
\begin{keywords}
Sparse representation, Mutual coherence, Gaussian random matrices, Bernoulli random matrices
\end{keywords}

\section{Introduction}\label{introduction}
\indent Recently, sparse representation has attracted a lot of scientists from many different fields. If $\mathbf{b}\in\mathbb{R}^{m}$ is a given signal which is known to be represented as a linear combination of only a few atoms of a dictionary $\mathbf{A}\in\mathbb{R}^{m\times n}$ (\textit{m}<\textit{n}), the representation vector $\mathbf{x}$ can be correctly and effectively computed by some procedures. The desire to find the maximally sparse solutions of an underdetermined linear system $\mathbf{Ax}=\mathbf{b}$ can be cast as the following optimization problem,
\begin{equation}\label{equation1}
(P_{0})\qquad\min_{\mathbf{x}}\,\|\mathbf{x}\|_{0}\quad s.t.\,\,\mathbf{b}=\mathbf{Ax}
\end{equation}

\section{Guarantees for Uniqueness and Stability} \label{guarantees}
\indent It is well known that that for a noiseless signal $\mathbf{b}_{0}=\mathbf{A}\mathbf{x}_{0}$, $\mathbf{x}_{0}$ is the unique solution to the problem (\ref{equation1}) with $\|\mathbf{x}\|_{0}<(1+1/\mu(\mathbf{A}))/2$ \cite{1} and both OMP \cite{2} and BP \cite{3} can find it \cite{1,4}. $\mu(\mathbf{A})$ is the largest absolute normalized inner product between different columns from $\mathbf{A}$. Clearly, the smaller the mutual coherence of $\mathbf{A}$ \cite{1} is, the better the above result becomes. \\
\indent An optimized algorithm was proposed in \cite{5} to reduce $\mu_{t}$ of the product of a projection and a given dictionary in the context of compressed sensing, a descendant of sparse approximation. \cite{6} presented improved versions of thresholding and OMP for sparse representation with an iterative algorithm which produced a new dictionary at the step of sensing with the cross cumulative coherence. While we work with the mutual coherence of an overcomplete Gaussian or Bernoulli random dictionary to improve the classical results in \cite{1,2,3,4}. \\
\indent Why we feel an interest in the mutual coherence of such dictionaries? Firstly, an overcomplete Gaussian or Bernoulli random matrix satisfies the RIP condition, another measure of a dictionary, which is widely used and different from the mutual coherence, with a high probability \cite{7}. So they have been commonly used in general numerical experiments to verify new methods. Secondly, the mutual coherence of such a matrix has been studied by some statistics from the limiting laws point of view \cite{8}. But how can a person directly reduces them in the finite dimensional case?

\section{Reducing the mutual coherence of an overcomplete Gaussian/Bernoulli random matrix} \label{reducing}
\indent When $\mathbf{A}$ is an overcomplete Gaussian or Bernoulli random matrix and has full row rank, the mutual coherence is fairly small (\ref{tail}) but can be reduced by multiplying an inverse matrix $\mathbf{P}$ generated from $\mathbf{A}$ itself on the left side. In this section, it is shown that the lower bound of the probability $\mathbb{P}\{\mu(\mathbf{PA})\leq\varepsilon\}$ with $\varepsilon\in(0,1)$ is larger than the one of $\mathbb{P}\{\mu(\mathbf{A})\leq\varepsilon\}$. The two bounds are obtained by the same way. \\
\indent A Bernoulli random matrix $\mathbf{A}$ means that $a_{ij}$ is independently selected to be $\pm$1/$\sqrt{m}$ with the same probability, while $a_{ij}$ of a Gaussian random matrix independently follows a normal distribution $N(0,\,m^{-1})$. $\bm\alpha_{j}$ denotes the \textit{j}\textsuperscript{th} column of $\mathbf{A}$. From the strong concentration of $\|\mathbf{Ax}\|_{2}^{2}$ in \cite{7}, it is known that the extreme singular values of $\mathbf{A}$, $\sigma_{1}$ and $\sigma_{m}$, satisfy
$0<\sqrt{1-\eta}\leq\sigma_{m}(\mathbf{A})\leq\sigma_{1}(\mathbf{A})\leq\sqrt{1+\eta}$ with the probability which is more than $1-2\exp(c_{0}(\eta))$ with $c_{0}(\eta)=\eta^{2}/4-\eta^{3}/6$, $\eta\in(0,1)$. So $\mathbf{A}$ has full rank with exponentially high probability. \\
\indent $\mathbf{P}=(\mathbf{A}\mathbf{A}^{T})^{-1/2}$ makes $(\mathbf{PA})(\mathbf{PA})^{T}=\mathbf{I}_{m}$. The intuitive idea to do this is that an orthogonal matrix has the smallest mutual coherence and $\mathbf{PA}$ consisting of \textit{m} orthonormal rows may have a smaller mutual coherence than $\mu(\mathbf{A})$. \\
\indent Consider $\mathbf{PA}$ and $\mathbf{A}$ in the singular value decomposition domain. The SVD of $\mathbf{A}$ is $\mathbf{A}=\mathbf{U}\,[\mathbf{S}, \mathbf{O}]\,\mathbf{V}^{T}$ with two orthonormal matrices $\mathbf{U}$ and $\mathbf{V}$ and a diagonal matrix, $\mathbf{S}=\operatorname{diag}\{\sigma_{1}\geq\sigma_{2}\geq\cdots\geq\sigma_{m}>0\}$. $\mathbf{V}$ can be divided into $\mathbf{V}_{1}\in\mathbb{R}^{n\times m}$ and $\mathbf{V}_{2}\in\mathbb{R}^{n\times(n-m)}$. $\mathbf{v}_{j}=[\mathbf{v}_{j}^{(1)},\mathbf{v}_{j}^{(2)}]$ is the \textit{j}\textsuperscript{th} row of $\mathbf{V}$ and $\mathbf{v}_{j}^{(1)}$ ,$\mathbf{v}_{j}^{(2)}$ are the \textit{j}\textsuperscript{th} rows of $\mathbf{V}_{1}$ and $\mathbf{V}_{2}$ respectively. The \textit{k}\textsuperscript{th} column of $\mathbf{A}=\mathbf{U}\mathbf{S}\mathbf{V}_{1}^{T}$ is $\bm\alpha_{k}=\mathbf{US}\mathbf{v}^{(1)T}_{k}$ and $\vert\langle\bm\alpha_{j},\bm\alpha_{k}\rangle\vert$ can be bounded by the Wielandt inequality in \cite{9,10} as
\begin{align*}
& \vert\langle\bm\alpha_{j},\bm\alpha_{k}\rangle\vert_{j\neq k}=\vert(\mathbf{US}\mathbf{v}_{j}^{(1)T})^{T}(\mathbf{US}\mathbf{v}_{k}^{(1)T})\vert=\vert\mathbf{v}^{(1)}_{j}\mathbf{S}^{2}\mathbf{v}^{(1)T}_{k}\vert \\ \leq & \left(\frac{\sigma_{1}^{2}-\sigma_{m}^{2}}{\sigma_{1}^{2}+\sigma_{m}^{2}}+\frac{\vert\langle\mathbf{v}^{(1)T}_{j},\mathbf{v}^{(1)T}_{k}\rangle\vert}{\|\mathbf{v}^{(1)}_{j}\|_{2}\|\mathbf{v}^{(1)}_{k}\|_{2}}\right)\cdot\left(1+\frac{\sigma_{1}^{2}-\sigma_{m}^{2}}{\sigma_{1}^{2}+\sigma_{m}^{2}}\cdot\frac{\vert\langle\mathbf{v}^{(1)T}_{j},\mathbf{v}^{(1)T}_{k}\rangle\vert}{\|\mathbf{v}^{(1)}_{j}\|_{2}\|\mathbf{v}^{(1)}_{k}\|_{2}}\right)^{-1}\triangleq u_{2}
\end{align*}
The \textit{k}\textsuperscript{th} column of $\mathbf{PA}$ is $\mathbf{P}\bm\alpha_{k}=\mathbf{U}\mathbf{v}^{(1)T}_{k}$ and $\vert\langle\mathbf{P}\bm\alpha_{j},\mathbf{P}\bm\alpha_{k}\rangle\vert/(\|\mathbf{P}\bm\alpha_{j}\|_{2}\|\mathbf{P}\bm\alpha_{k}\|_{2})$ which is equal to $\vert\langle\mathbf{v}^{(1)T}_{j},\mathbf{v}^{(1)T}_{k}\rangle\vert/(\|\mathbf{v}^{(1)}_{j}\|_{2}\|\mathbf{v}^{(1)}_{k}\|_{2})\triangleq u_{3}$.
Clearly, $u_{2}$ is larger than $u_{3}$. For the subset of overcomplete Gaussian and Bernoulli random matrices with the aforementioned positive bounded singular values, consider the following events $E_{1}=\{(1-\eta)\|\mathbf{x}\|_{2}^{2}\leq\|\mathbf{Ax}\|_{2}^{2}\leq(1+\eta)\|\mathbf{x}\|_{2}^{2},\exists\eta\in(0,1)\}$, $E_{2}=\{\mu(\mathbf{A})\leq\varepsilon,\varepsilon\in(0,1)\}$ and $E_{3}=\{\mu(\mathbf{PA})\leq\varepsilon,\varepsilon\in(0,1)\}$. It is shown that
\begin{align*}
&\mathbb{P}\{E_{1}\cap E_{2}\}=\mathbb{P}\{E_{1}\}-\mathbb{P}\{E_{1}\cap E_{2}^{c}\}\quad
\mathbb{P}\{E_{1}\cap E_{3}\}=\mathbb{P}\{E_{1}\}-\mathbb{P}\{E_{1}\cap E_{3}^{c}\}\\
&\mathbb{P}\{E_{1}\cap E_{2}^{c}\}\leq n(n-1)\mathbb{P}\{\varepsilon\leq u_{2}\}\triangleq p_{2}\quad
\mathbb{P}\{E_{1}\cap E_{3}^{c}\}\leq n(n-1)\mathbb{P}\{\varepsilon\leq u_{3}\}\triangleq p_{3}
\end{align*}
\indent So it indirectly reflects the phenomenon that the mutual coherence of the product of an inverse matrix $\mathbf{P}$ and an overcomplete Gaussian or Bernoulli random matrix is smaller than the original mutual coherence from the result that the former has a better lower bound of the probability $\mathbb{P}\{\mu(\mathbf{PA})\leq\varepsilon\}$ than the one of $\mathbb{P}\{\mu(\mathbf{A})\leq\varepsilon\}$, $\varepsilon\in(0,1)$, although they are obtained by the same way.

\section{The essential mutual coherence} \label{essential}
In this section, we use the above way for reducing the mutual coherence of an overcomplete Gaussian or Bernoulli dictionary on a set $\mathbf{X}$ which almost includes the two types and prove that the newly defined essential mutual coherence on $\mathbf{X}$ is strictly smaller than the original one but a Lebesgue zero measure subset. \\
\indent $\mathbf{X}\triangleq\{\mathbf{A}\in\mathbb{R}^{m\times n};m<n,\mu(\mathbf{A})<1,0\not\in\operatorname{diag}(\mathbf{A}^{T}\mathbf{A})\}$. \textit{m}\,<\,\textit{n} is natural and $\mathbf{A}$ has no zero columns otherwise there exist some meaningless atoms. It is senseless to multiply an inverse matrix $\mathbf{P}$ from the left side of $\mathbf{A}$ for reducing $\mu(\mathbf{A})$, if $\mu(\mathbf{A})=1$. For the Gaussian type, $\mathbb{P}\{0\in\operatorname{diag}(\mathbf{A}^{T}\mathbf{A})\}=0=\mathbb{P}\{\mu(\mathbf{A})=1\}$ if we notice that the independence of each $a_{ij}$ and the condition on equality of Cauchy Schwarz Inequality. For the Bernoulli type, $\mathbb{P}\{\mu(\mathbf{A})=1\}$\,$\leq$\,$n(n-1)\exp(-m/2)$ due to the Hoeffding Inequality in probability. That is to say that the Gaussian or Bernoulli type belongs to $\mathbf{X}$ with a high probability.\\
\indent Consider an equivalent problem of $(P_{0})$, $(P_{0}^{\prime})$\, $\min_{\mathbf{\mathbf{x}}}\,\|\mathbf{x}\|_{0}\,s.t.\,\mathbf{Pb}=\mathbf{PAx}$ with an invertible matrix $\mathbf{P}$. A new quantity, essential mutual coherence, is defined as $\mu_{e}(\mathbf{A})\triangleq\inf_{\mathbf{P}}\mu(\mathbf{PA})$, $\mathbf{P}\in GL_{m}(\mathbb{R})=\{\text{all \textit{m} order real invertible matrices}\}$ and is invariant for the elementary row operations of matrices. In this section we show that $\mu_{e}(\mathbf{A})<\mu(\mathbf{A})$ holds true almost every where on $\mathbf{X}$ by constructing a matrix $\mathbf{P}\triangleq\mathbf{E}_{m}+\varepsilon\mathbf{E}_{m,1}\in GL_{m}(\mathbb{R})$ with a proper $\varepsilon$ such that $\mu(\mathbf{PA})<\mu(\mathbf{A})$. $\mathbf{E}_{m}$ is the \textit{m} order identity matrix and $\mathbf{E}_{m,1}\in\mathbb{R}^{m\times m}$ consists of 0 but 1 on the position (\textit{m}, 1). \\
\indent How to choose the parameter $\varepsilon$? For two different columns of $\mathbf{A}$, $\bm\alpha_{i}$ and $\bm\alpha_{j}$, $i<j$, a polynomial $f_{ij}(\varepsilon)=A_{ij}+B_{ij}\varepsilon+C_{ij}\varepsilon^{2}+D_{ij}\varepsilon^{3}+E_{ij}\varepsilon^{4}$ can be obtained from $I(\mathbf{P}\bm\alpha_{i},\mathbf{P}\bm\alpha_{j})\triangleq\vert\langle\mathbf{P}\bm\alpha_{i},\mathbf{P}\bm\alpha_{j}\rangle\vert/(\|\mathbf{P}\bm\alpha_{i}\|_{2}\|\mathbf{P}\bm\alpha_{j}\|_{2})<\mu(\mathbf{A})$ (use $\mu$ instead of $\mu(\mathbf{A})$).\\
\indent $f_{ij}(0)=A_{ij}=\|\bm\alpha_{i}\|_{2}^{2}\|\bm\alpha_{j}\|_{2}^{2}(I(\bm\alpha_{i},\bm\alpha_{j})^{2}-\mu^{2})<0$, for all pairs ($\bm\alpha_{i}$,\,$\bm\alpha_{j}$) such that $I(\bm\alpha_{i},\bm\alpha_{j})<\mu$. There exists a $\varepsilon^{(ij)}_{1}>0$ such that for all $\varepsilon\in(-\varepsilon^{(ij)}_{1},\,\varepsilon^{(ij)}_{1})$, $f_{ij}(\varepsilon)<0$ owing to continuity of the polynomial $f_{ij}(\varepsilon)$. A minimum, $\varepsilon_{1}$, can be founded among all the $\varepsilon^{(ij)}_{1}>0$ for all such pairs. \\
\indent For all pairs $(\bm\alpha_{i}$,\,$\bm\alpha_{j})$ such that $I(\bm\alpha_{i},\bm\alpha_{j})=\mu$, if $f_{ij}^{'}(0)=B_{ij}$ is positive or negative, then there exists a $\varepsilon_{2}^{(ij)}>0$ such that $f_{ij}(\varepsilon)<0$, $\forall\varepsilon\in(-\varepsilon_{2}^{(ij)},0)$ or $(0,\varepsilon^{(ij)}_{2})$. Among all of them there exists a minimum, $\varepsilon_{2}$. So the sign of the final $\varepsilon$ is attributed to all of the pairs such that $I(\bm\alpha_{i},\bm\alpha_{j})=\mu$.\\
\indent Two exceptions catch our attention. Firstly, it is impossible to choose $\varepsilon_{1}^{(ij)}$ or $\varepsilon_{2}^{(ij)}$ when both $A_{ij}$ and $B_{ij}$ are zero. The set $\mathbf{X}_{1}\triangleq\{\mathbf{A}\in\mathbf{X};B_{ij}=0, A_{ij}=0\}$ has measure zero due to $\mathbf{X}_{1}\subset\cup_{i<j}\mathbf{E}_{ij}$ and $\mathbf{E}_{ij}$ which consists of all matrices that satisfy $\|\bm\alpha_{i}\|_{2}^{2}\|\bm\alpha_{j}\|_{2}^{2}(a_{mi}a_{1j}+a_{1i}a_{mj})-\langle\bm\alpha_{i},\bm\alpha_{j}\rangle(\|\bm\alpha_{i}\|_{2}^{2}a_{1j}a_{mj}+\|\bm\alpha_{j}\|_{2}^{2}a_{1i}a_{mi})=0$ is a Lebesgue measure zero set (\ref{zeromeasure}), if matrix $\mathbf{A}\in\mathbf{X}$ is seen as a vector in $\mathbb{R}^{mn}$ in some order. So $A_{ij}=0$ or $B_{ij}=0$ seldom occur at the same time.\\
\indent Secondly, $\varepsilon_{2}$ can not be obtained if there exist two different pairs such that $I(\bm\alpha_{i_{1}},\bm\alpha_{j_{1}})=I(\bm\alpha_{i_{2}},\bm\alpha_{j_{2}})=\mu$ and $B_{i_{1}j_{1}}B_{i_{2}j_{2}}<0$. Consider a set that contains more elements, $\mathbf{X}_{2}\triangleq\{\mathbf{A}\in\mathbf{X};\mu=I(\bm\alpha_{i_{k}},\bm\alpha_{j_{k}})_{i_{k}<j_{k}},k>1\}$, a subset of $\cup_{(i_{1},j_{1})\neq(i_{2},j_{2})}\mathbf{F}_{(i_{1},j_{1}),(i_{2},j_{2})}$. Similarly, $\mathbf{F}_{(i_{1},j_{1}),(i_{2},j_{2})}=\{\mathbf{A},\vert\langle\bm\alpha_{i_{1}},\bm\alpha_{j_{1}}\rangle\vert^{2}\|\bm\alpha_{i_{2}}\|_{2}^{2}\|\bm\alpha_{j_{2}}\|_{2}^{2}=\vert\langle\bm\alpha_{i_{2}},\bm\alpha_{j_{2}}\rangle\vert^{2}\|\bm\alpha_{i_{1}}\|_{2}^{2}\|\bm\alpha_{j_{1}}\|_{2}^{2}\}$ is also a measure zero set. So only one pair almost always achieve $\mu$.\\
\indent When $\mathbf{A}\in\mathbf{X}$ is given, an interval $(-\varepsilon^{(ij)}_{1},\varepsilon^{(ij)}_{1})$ can be selected for the pairs such that $I(\bm\alpha_{i},\bm\alpha_{j})<\mu$ and an interval $(-\varepsilon_{2}^{(ij)},0)$ or $(0,\varepsilon_{2}^{(ij)})$ can be chosen if $\bm\alpha_{i}$ and $\bm\alpha_{j}$ satisfy $I(\bm\alpha_{i},\bm\alpha_{j})=\mu$. So a proper $\varepsilon$ which is in the intersection among all of the above intervals can be selected to obtain $\mu(\mathbf{PA})<\mu(\mathbf{A})$. All the exceptions that are included in $\mathbf{X}$ fall into $\mathbf{X}_{1}\cup\mathbf{X}_{2}$, which has measure zero. So a proper $\varepsilon$ with plus or minus sign can be chosen to construct $\mathbf{P}$ such that $\mu(\mathbf{PA})<\mu(\mathbf{A})$ almost everywhere on the set $\mathbf{X}$.

\section{Experiments and results} \label{experiment}
\begin{figure}[htbp]
\centering
\includegraphics[width=9cm,height=5cm]{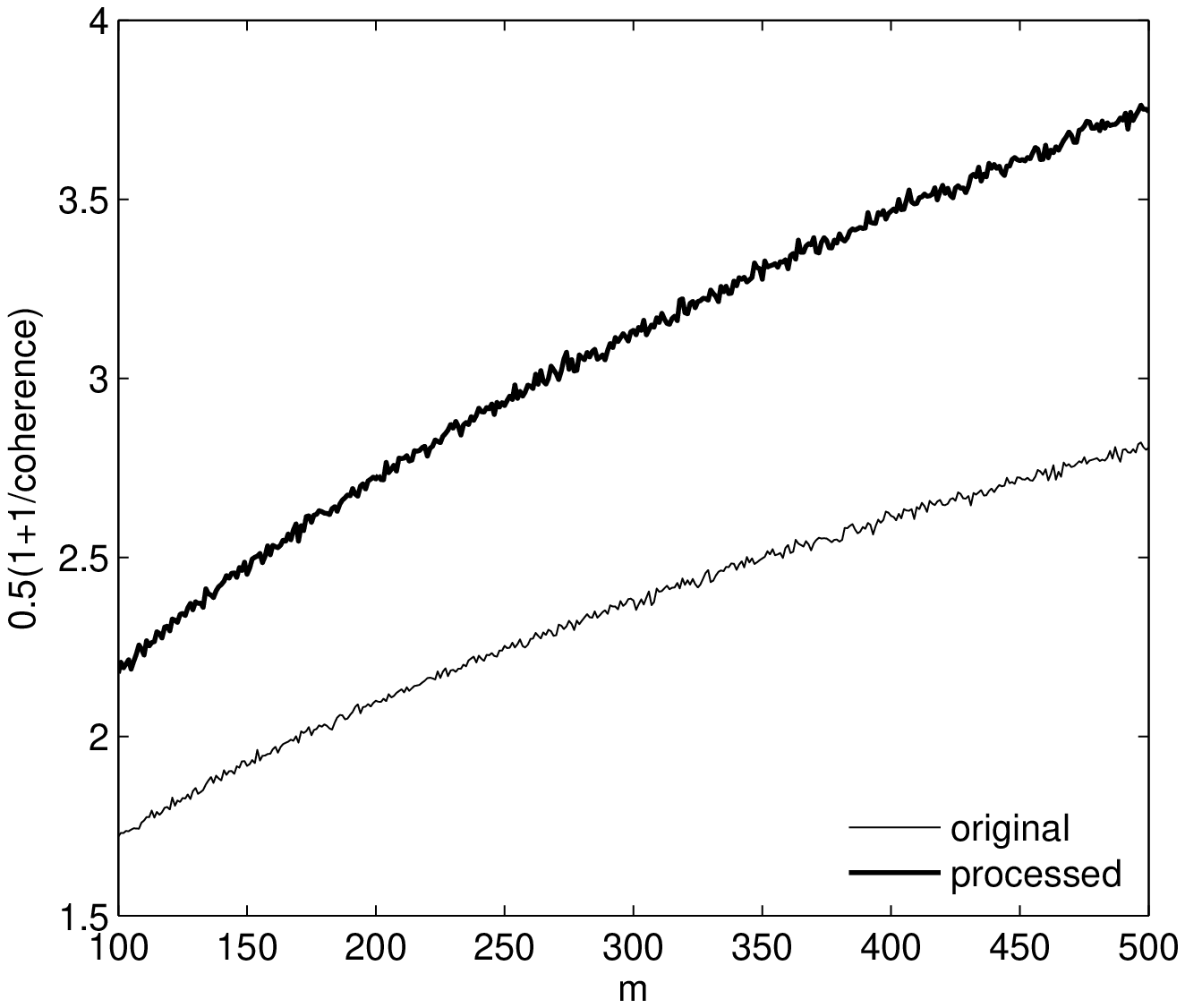}
\captionsetup{aboveskip=0pt,belowskip=0pt,labelfont=bf,labelsep=period}
\renewcommand{\figurename}{Fig.}
\caption{A comparison of $0.5(1+1/\mu(\mathbf{A}))$ and $0.5(1+1/\mu(\mathbf{PA}))$ for the Bernoulli random dictionary of size $m\times 2m$, where $m\in[100, 500]$.}
\label{figure1}
\end{figure}

\begin{figure}[htbp]
\centering
\includegraphics[width=9cm,height=5cm]{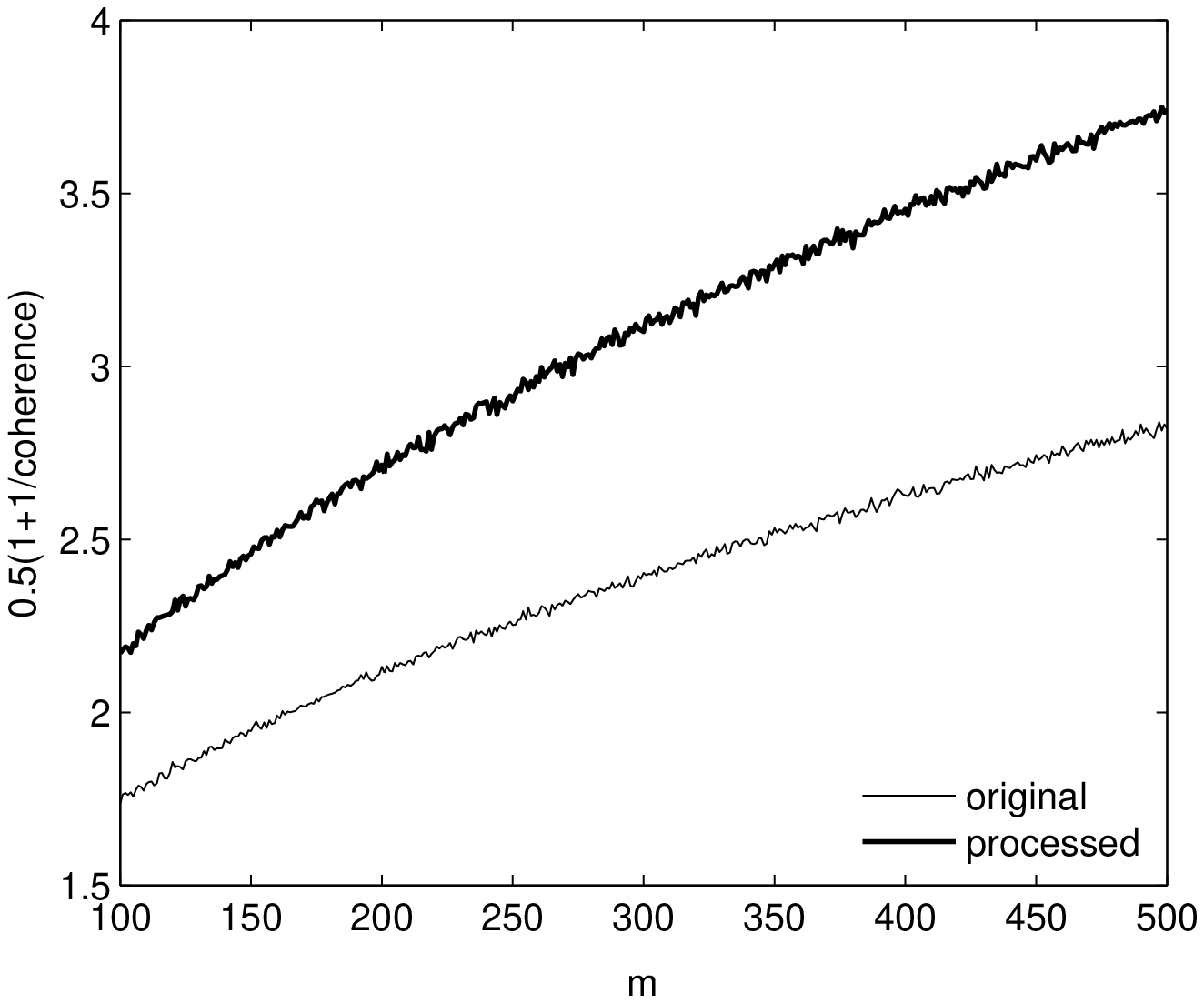}
\captionsetup{aboveskip=0pt,belowskip=0pt,labelfont=bf,labelsep=period}
\renewcommand{\figurename}{Fig.}
\caption{A comparison of $0.5(1+1/\mu(\mathbf{A}))$ and $0.5(1+1/\mu(\mathbf{PA}))$ for the Gaussian random dictionary of size $m\times 2m$, where $m\in[100, 500]$.}
\label{figure2}
\end{figure}

\begin{figure}[htbp]
\centering
\includegraphics[width=9cm,height=5cm]{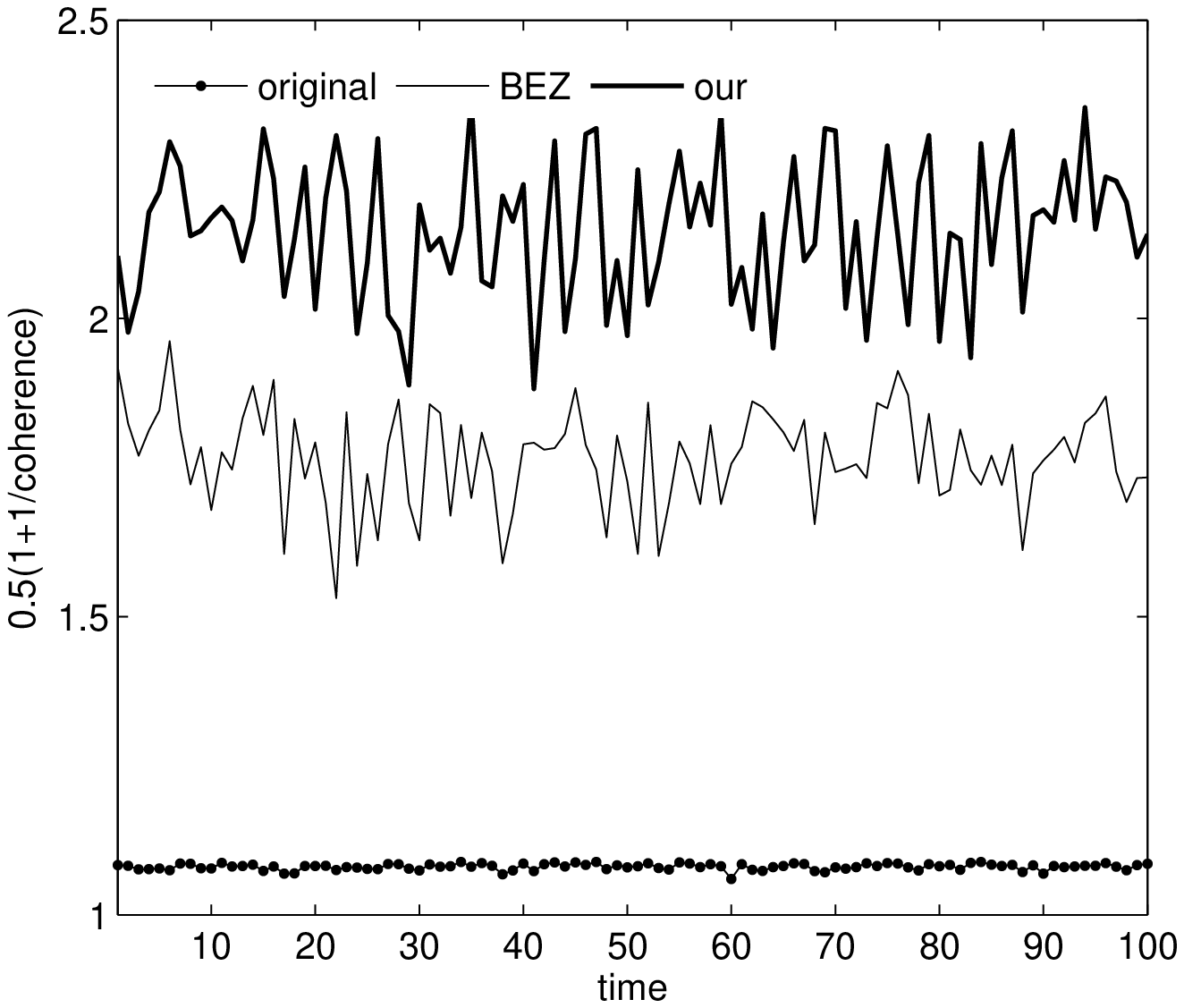}
\captionsetup{aboveskip=0pt,belowskip=0pt,labelfont=bf,labelsep=period}
\renewcommand{\figurename}{Fig.}
\caption{A comparison of 0.5(1+1/$\mu(\mathbf{D})$), 0.5(1+1/$\mu(\mathbf{PD})$) with $\mathbf{P}$\,=\,($\mathbf{D}\mathbf{D}^{T})^{-1/2}$ and 0.5(1+1/$\mu(\mathbf{PD})$) with $\mathbf{P}$\,=\,$(1-\epsilon)/100\mathbf{1}\mathbf{1}^{T}$ for $\mathbf{D}$ of size $100\times200$.}
\label{figure3}
\end{figure}

\indent In the first experiment, both $\mu(\mathbf{A})$ and $\mu(\mathbf{PA})$ are calculated for the $m\times2m$ Bernoulli (Fig. \ref{figure1}) or Gaussian (Fig. \ref{figure2}) random dictionary averaged over 100 times for each value of \textit{m} in the range [100, 500].\\
\indent Only Fig. \ref{figure1} is referred owing to the similarity between the two figures. From about \textit{m}\,=\,280, the bound 0.5(1+1/coherence) exceeds 3 and keeps increasing till 500 when coherence is taken as $\mu(\mathbf{PA})$ but the one with $\mu(\mathbf{A})$ remains below 3 till 500 although it also becomes larger and larger. That is, $\mu(\mathbf{PA})$ can lead to the recovery of the sparse solution with 3 nonzero elements for both OMP and BP when $m\in[280, 500]$. It is one more than the result obtained by using $\mu(\mathbf{A})$.
Both $0.5(1+1/\mu(\mathbf{A}))$ and $0.5(1+1/\mu(\mathbf{PA}))$ grow alongside the increase in \textit{m}, so only six values of \textit{m} are listed in Table \ref{table1} in order to avoid the similarity between the forms of Fig. \ref{figure1} and Table \ref{table1} which considers all the values in [1500,\,2000]. As shown in Table \ref{table1}, the condition $0.5(1+/\mu(\mathbf{A}))$ is below 5 but $0.5(1+1/\mu(\mathbf{PA}))$ is above 6 when \textit{m} is after 1800.
\begin{table}[htbp]
\captionsetup{labelfont=bf,font=normalfont,labelsep=newline,justification=RaggedRight,aboveskip=0pt,belowskip=0pt}
\caption{A comparison of 0.5(1+1/$\mu(\mathbf{A})$) and 0.5(1+1/$\mu(\mathbf{PA})$) for the Bernoulli random matrix of size $m\times 2m$ when $m\in\operatorname{linspace}(1500, 2000, 100)$ (a MATLAB$^{\circledR}$ notation).}\label{table1}
\arrayrulewidth=1pt
\tabcolsep=7.5pt
\begin{flushleft}
  \begin{tabular}{l l l l l l l}
  \hline
   %after \\: \hline or \cline{col1-col2} \cline{col3-col4} ...
   \textit{m}         & 1500   & 1600   & 1700   & 1800   & 1900   & 2000 \\ \hline
   0.5(1+1/$\mu(\mathbf{A})$)  & 4.1531 & 4.2708 & 4.3896 & 4.4893 & 4.5769 & 4.6738 \\
   0.5(1+1/$\mu(\mathbf{PA})$) & 5.6711 & 5.8311 & 5.9798 & 6.0945 & 6.2611 & 6.4046 \\
  \hline
\end{tabular}
\end{flushleft}
\end{table}
\begin{table}[htbp]
\captionsetup{labelfont=bf,font=normalfont,labelsep=newline,justification=RaggedRight,aboveskip=0pt,belowskip=0pt}
\caption{A comparison of 0.5(1+1/$\mu(\mathbf{A})$) and 0.5(1+1/$\mu(\mathbf{PA})$) for the Gaussian random matrix of size $m\times 2m$ when $m\in\operatorname{linspace}(1500, 2000, 100)$.}\label{table2}
\arrayrulewidth=1pt
\tabcolsep=7.5pt
\begin{flushleft}
  \begin{tabular}{l l l l l l l}
  \hline
   %after \\: \hline or \cline{col1-col2} \cline{col3-col4} ...
   \textit{m}         & 1500   & 1600   & 1700   & 1800   & 1900   & 2000 \\ \hline
   0.5(1+1/$\mu(\mathbf{A})$)  & 4.1686 & 4.3185 & 4.4032 & 4.4860 & 4.5952 & 4.6469 \\
   0.5(1+1/$\mu(\mathbf{PA})$) & 5.7126 & 5.8191 & 5.9664 & 6.0822 & 6.2473 & 6.4228 \\
  \hline
\end{tabular}
\end{flushleft}
\end{table}
\indent The second experiment compares the effect of the above method and the one proposed in \cite{11} (called BEZ here) for decreasing the mutual coherence of the overcomplete uniform random dictionary $\mathbf{D}$ which is obtained by choosing the entries independently from a uniform distribution in [0,\,1] and then normalizing the columns to a unit $\ell_{1}$-norm. The authors in \cite{11} used $\mathbf{P}=(1-\epsilon)/100\mathbf{1}\mathbf{1}^{T}$. Here, $\mathbf{1}\in\mathbb{R}^{100}$ is of 1's and $\epsilon$ satisfies $0<\epsilon\ll1$. \\
\indent Three values, $0.5(1+1/\mu(\mathbf{D}))$, $0.5(1+1/\mu(\mathbf{PD}))$ with $\mathbf{P}=(\mathbf{D}\mathbf{D}^{T})^{-1/2}$ and $0.5(1+1/\mu(\mathbf{PD})$) with $\mathbf{P}=(1-\epsilon)/100\mathbf{1}\mathbf{1}^{T}$ are obtained for the dictionary $\mathbf{D}$ generated as per instructions in \cite{11} with a fixed size of $100\times200$. This is done 100 times and $\epsilon$ is selected to maximize $0.5(1+1/\mu(\mathbf{PD})$) with $\mathbf{P}$ in \cite{11} over $\operatorname{linspace}(10^{-4},10^{-1},1000)$ per time. In Fig. \ref{figure3}, three lines placed from bottom to up represent the aforementioned three values respectively. $0.5(1+1/\mu(\mathbf{D}))$ is below 1.5 all the time. In all 100 times, $0.5(1+1/\mu(\mathbf{PD}))$ obtained by using BEZ in \cite{11} is strictly smaller than 2, but meanwhile, the one generated by applying our method can leap two 86 times.

\section{Conclusion} \label{conclusion}
\indent The mutual coherence of an overcomplete Gaussian or Bernoulli random matrix which is proven to be small can be reduced directly by multiplying an inverse matrix from the left side to improve the traditional umbrella of unique sparse representation and successful reconstruction behavior. Furthermore, the newly proposed essential mutual coherence which is inspired by what has been done is proven to be strictly smaller than the original one on a set except a Lebesgue measure zero set. Numerical results exhibits the decrease in the mutual coherence of an overcomplete Gaussian or Bernoulli random dictionary and also an overcomplete Uniform random matrix  which is better than the former result obtained by other authors.

\begin{appendix}
\section{Tail bound}\label{tail}
\indent For Gaussian random matrices, $m\|\bm\alpha_{j}\|_{2}^{2}\sim\chi^{2}(m)$.
\begin{itemize}
  \item $\mathbb{P}\{Z-m\leq -2\sqrt{mx}\}\leq\exp(-x)$, $\forall x>0$, for a centralized $\chi^{2}$-variable \textit{Z} with \textit{m} degrees of freedom \cite{12}.
  \item $\mathbb{P}\{\vert\langle\bm\alpha_{j},\bm\alpha_{k}\rangle\vert_{j<k}\geq x\}\leq2\exp(-0.25mx^{2}/(1+0.5x))$, $\forall x>0$ owing to the independence of $\bm\alpha_{j}$ and $\bm\alpha_{k}$, $j<k$ \cite{13}.
\end{itemize}
For all $a\in(0,1)$, according to the fact that
$$\mathbb{P}\{\mu(\mathbf{A})\geq\varepsilon\}\leq n(n-1)/2\left(\mathbb{P}\{\vert\langle\bm\alpha_{j},\bm\alpha_{k}\rangle\vert_{j<k}\geq a\varepsilon\}+2\mathbb{P}\{\|\bm\alpha_{j}\|_{2}^{2}\leq a\}\right)$$ we have
\[
\mathbb{P}\{\mu(\mathbf{A})\geq\varepsilon\in(0,1)\}\leq n(n-1)/2\left[\exp\left(-\frac{ma^{2}\varepsilon^{2}}{4(1+a\varepsilon/2)}\right)+\exp\left(-\frac{m}{4}(1-a)^{2}\right)\right]
\]
\indent A similar tail bound can be obtained by applying the Hoeffding Inequality in probability for the Bernoulli case

\section{$\mathbf{E}_{ij}$ is a Lebesgue measure zero set}\label{zeromeasure}
\indent Without loss of generality, we show that $m(E_{12})=0$ (\textit{m} denotes the Lebesgue measure on $\mathbb{R}^{mn}$). $\mathbf{E}_{12}\subset\mathbf{R}_{12}\times\mathbb{R}^{m\times(n-2)}\triangleq\mathbf{S}_{12}=\cup_{r}^{\infty}\mathbf{S}_{r}$. $\mathbf{R}_{12}$ is a set which consists of all $(\bm\alpha_{1}^{T},\bm\alpha_{2}^{T})^{T}$ that satisfy
\begin{center}
  $\|\bm\alpha_{1}\|_{2}^{2}\|\bm\alpha_{2}\|_{2}^{2}(a_{m1}a_{12}+a_{11}a_{m2})-\langle\bm\alpha_{1},\bm\alpha_{2}\rangle(\|\bm\alpha_{1}\|_{2}^{2}a_{12}a_{m2}+\|\bm\alpha_{2}\|_{2}^{2}a_{11}a_{m1})=0$
\end{center}
$\mathbf{S}_{r}=\mathbf{S}_{r}^{(1)}\times\mathbf{S}_{r}^{(2)}$ with $\mathbf{S}_{r}^{(2)}=(-r,r)^{m(n-2)}$, $r=1,2,3,\cdots$ and $\mathbf{S}_{r}^{(1)}=\mathbf{R}_{12}$. It is known that $m(\mathbf{S}_{r})=m(\mathbf{S}_{r}^{(1)})\cdot m(\mathbf{S}_{r}^{(2)})=0$ due to $m(\mathbf{S}^{(1)}_{r})=0$ and $m(\mathbf{S}^{(2)}_{r})=(2r)^{m(n-2)}$. So $m(\mathbf{S}_{12})=\lim_{r\rightarrow\infty}m(\mathbf{S}_{r})=0$ because $\mathbf{S}_{1}\subset\mathbf{S}_{2}\subset\cdots\subset\mathbf{S}_{r}\subset\cdots$.\\
\indent $m(\mathbf{R}_{12})=0$ due to a polynomial from $\mathbb{R}^{n}$ to $\mathbb{R}$ is either identically zero or nonzero almost everywhere.
\end{appendix}

\bibliographystyle{plain}
\bibliography{manuscript}

\end{document}